\documentclass[12pt]{amsart}
\usepackage[utf8]{inputenc}
\usepackage[T2A]{fontenc}
\usepackage{textcomp}
\usepackage{amsmath, amssymb}

\usepackage[left=3cm,right=3cm,
    top=3cm,bottom=3cm,bindingoffset=0cm]{geometry}
 
\begin{document}

\title{Rota-Baxter operators of non-scalar weights, connections with coboundary Lie bialgebra structures.}
\author{{M. Goncharov}}%
\address{Maxim Goncharov
\newline\hphantom{iii} Novosibirsk State University, Department of Mechanics and
Mathematics, Novosibirsk, Russia}%
\email{goncharov.gme@gmail.com}%


\vspace{1cm}
\maketitle {\small
\begin{quote}
\noindent{\sc Abstract. } In the paper, we introduce the notion of a Rota-Baxter operator of a non-scalar weight. As a motivation, we show that there is a natural connection between Rota-Baxter operators of this type and structures of quasitriangular Lie bialgebras on a quadratic finite-dimensional Lie algebra. Moreover, we show that some classical results on Lie bialgebra structure on simple finite-dimensional Lie algebras can be obtained from the corresponding results for Rota-Baxter operators.\medskip

\noindent{\bf Keywords:} Lie algebra, Lie bialgebra, coboundary Lie bialgebra, classical Yang-Baxter equation, Rota---Baxter operator
 \end{quote}
}

Rota---Baxter operators for associative and commutative algebras first appeared in the
paper by G. Baxter \cite{Baxter} as a tool for studying integral operators that
appear in the theory of probability and mathematical statistics. The original notion of the Rota-Baxter operators was the following: given an associative and commutative algebra $A$, by a Baxter type operator one meant a linear map $R:A\rightarrow A$ satisfying the following quadratic equation
$$
R(a)R(b)=R(R(a)b+aR(b)+\theta ab),
$$
where $\theta\in A$ is a fixed element, $a,b\in A$ \cite{Atkinson}. For a long period of time, this type of operators had been intensively studied in combinatorics and probability theory, mainly.

Independently, Rota---Baxter  operators on Lie algebras
naturally appeared in the papers of A.A. Belavin, V.G. Drinfeld \cite{BD}, M.A. Semenov-Tyan-Shanskii \cite{STS} and N. Yu. Reshetikhin, M.A. Semenov-Tyan-Shanskii \cite{STSRESH} while studying triangular and factorizable Lie bialgebras structures on simple finite-dimensional complex Lie algebras (see Section 3 for details). In this context, the definition of a Rota-Baxter operator was the following: given an arbitrary algebra $A$ over a field $F$ and $\lambda\in F$, by a Rota-Baxter operator of weight $\lambda$ one meant a linear map $R:A\rightarrow A$ satisfying
\begin{equation}\label{RBold}
R(a)R(b)=R(R(a)b+aR(b)+\lambda ab)
\end{equation}
for all $a,b\in A$. Nowadays, Rota-Baxter operators of weight $\lambda$ are defined as operators satisfying \eqref{RBold}. Operators of this type turn out to be interesting, with many interesting applications in various areas of mathematics, such as symmetric polynomials, quantum field
renormalization, pre- and postalgebras, shuffle algebra, double Poisson algebras,  etc. (see \cite{Aguiar}, \cite{Atkinson}, \cite{Fard}, \cite{Guo}, \cite{Ogievetsky}, \cite{GonGub}).

At the same time, the notion of a Rota-Baxter operator for groups was introduced in \cite{GLY}. If $G$ is a group, then a map $B:G\mapsto G$ is called a Rota-Baxter operator on the group $G$ if for all $g,h\in G$:
$$
B(g)B(h)=B(gB(g)hB(g)^{-1}).
$$
If $(G,B)$ is a Rota-Baxter Lie group, then the tangent map of $B$ at the identity is a Rota-Baxter operator of weight 1 on the corresponding Lie algebra of the Lie group $G$. Also, it was shown that many results that are true for Rota-Baxter operators on algebras have corresponding analogs for Rota-Baxter operators on groups. And the natural question here is: what is a weight in the case of Rota-Baxter operators on groups, and how to define a Rota-Baxter operator of any weight? Considering elements of a field as a weight in this case is not natural. The same question can be asked for other algebraic structures ( rings, Hopf algebras, etc.).

In this paper, we consider a generalization of the notion of a weight of a Rota-Baxter of an algebra. The main motivation of the paper is to provide a definition of a Rota-Baxter operator that can be used for the study of quasitriangular (not necessarily factorizable) Lie bialgebras. Instead of elements from a ground field $F$ (or elements of an associative and commutative algebra), we suggest considering centrum maps, that is, maps that commute with operators of left and right multiplications (see Section 2). This is a generalization of both notions of Rota-Baxter operators: given a scalar $\lambda\in F$ (or an element $\theta\in A$ if $A$ is an associative commutative algebra), we can consider a linear map $\varphi_{\lambda}$ (respectively, $\varphi_{\theta})$ acting as $\varphi_{\lambda}(x)=\lambda x$. From this point of view, the Rota-Baxter operator on groups from \cite{GLY} can be called Rota-Baxter operators of weight $\mathrm{id}$ (where $\mathrm{id}$ is the identity map).

The paper is organized as follows: In Section 1, we give preliminary results concerning Lie bialgebras and the classical Yang-Baxter equation on Lie algebras. In Section 2, the notion of a Rota-Baxter operator of a non-scalar weight $\mu\in \mathrm{Cent}$ is introduced and studied. In Section 3, we obtain the connection between Rota-Baxter operators on a finite-dimensional Lie algebra $\mathfrak g$ over an arbitrary field of characteristic different from 2 and solutions of the classical Yang-Baxter equations on $\mathfrak g$ with $\mathfrak g$-invariant symmetric part (that is, with structures of triangular and quasitriangular Lie bialgebra structures on $\mathfrak g$). In Section 4, we use the results from Section 2 to study coboundary Lie bialgebra structures on perfect finite-dimensional quadratic Lie algebras over an arbitrary field of characteristic different from 2. Given $r\in \mathfrak g\wedge\mathfrak g$, we establish a criteria for a coalgebra $(\mathfrak g, \delta_r)$ to be a Lie bialgebra. Moreover, we can use the obtained results to specify whether $(\mathfrak g, \delta_r)$ is a triangular, quasitriangular, or factorizable Lie bialgebra. Finally, in Section 5, we show that many well-known results concerning Lie bialgebra structures on simple finite-dimensional Lie bialgebras over a field of characteristic 0 can be obtained as corollaries of the results from Sections 2-4. Finally, we show that every Lie bialgebra structure on a semisimple finite-dimensional Lie bialgebra over a field of characteristic 0 is induced by a special type of Rota-Baxter operator of some (not necessarily a  scalar or even a non-degenerate) weight.

The author acknowledges the kind hospitality of Jilin University. The author is also grateful to professors Y. Sheng and R. Tang for useful discussions.

\section{Preliminaries: Lie bialgebras and CYBE.}
Unless otherwise specified, $F$ denotes a field of characteristic different from 2. Given vector spaces $V$ and $U$ over a field $F$, denote by $V\otimes U$ its tensor product over $F$. Define the linear mapping $\tau$ on $V\otimes V$ by $\tau(\sum\limits_ia_i\otimes
b_i)=\sum\limits_ib_i\otimes a_i$. If $\tau(r)=-r$, then an element $r\in V\otimes V$ is called skew-symmetric, and we will use the notation $V\wedge V$ for the subspace of skew-symmetric elements in $V\otimes V$.

Denote by $V^*$ the dual space of $V$. Given $f\in V^*$ and $v\in V$, the symbol $\langle f,v\rangle$   (or $\langle v,f\rangle$ )denotes the the linear functional $f$ evaluated on $v$ (i.e. $\langle f,v\rangle=\langle v,f\rangle =f(v)$). 

Given a Lie algebra $\mathfrak g$, denote by $Z(\mathfrak g)$ the center of the algebra $\mathfrak g$:
$$
Z(\mathfrak g)=\{x\in \mathfrak g|\ [x,\mathfrak g]=0\}.
$$ 

{\bf Definition.} A lie algebra $\mathfrak g$ is called perfect, if $[\mathfrak g,\mathfrak g]=\mathfrak g$.

A Lie algebra $\mathfrak g$ acts on $\mathfrak g^{\otimes n}$ by
$$
[x_1\otimes x_2\otimes \ldots \otimes x_n,y]=\sum\limits_i
x_1\otimes\ldots\otimes [x_i,y]\otimes \ldots\otimes x_n$$

for all $x_i,y\in \mathfrak g$.

An element $r\in L^{\otimes n}$ is called $\mathfrak g$-invariant (or $ad$-invariant) if $[r,y]=0$ for all $y\in L$.

{\bf Definition.} Let $A$ be an arbitrary vector space over $F$. A pair $(A, \delta)$, where 
$\delta : A \rightarrow A\otimes A$ is a linear mapping, is called a
coalgebra, while $\delta$ is a comultiplication.

Given $a\in A$, we will use the following Sweedler's sumless notation: 
$$\delta(a) =  a_{(1)}\otimes a_{(2)}.$$

Let $(A,\delta)$ be a coalgebra. Define a multiplication on $A^*$ by
$$
 fg(a)=\langle fg,a\rangle= f(a_{(1)})g(a_{(2)}),
$$
where $f,g\in A^*$, $a\in A$ and $\delta(a)= a_{(1)}\otimes
a_{(2)}$.  The algebra obtained is  \emph{the dual algebra} of the
coalgebra $(A, \delta)$.

The following definition of a coalgebra related to some variety of algebras was given in \cite{ANQ}.

{\bf Definition.} Let $\mathcal{M}$ be an arbitrary variety of
algebras. The pair $(A, \delta)$ is called a $\mathcal{M}$-coalgebra
if $A^*$ belongs to $\mathcal{M}$.

The following definition was given in \cite{Drinf}:

{\bf Definition} Let $\mathfrak g$ be a Lie algebra with a comultiplication
$\delta$. The pair $(\mathfrak g,\delta)$ is called a Lie bialgebra if and
only if $(\mathfrak g,\delta)$ is a Lie coalgebra and $\delta$ is a 1-cocycle,
i.e., it satisfies
\begin{equation}\label{bder}
\delta([a,b])=[\delta(a),b]+[a,\delta(b)]
\end{equation}
for all $a,b\in \mathfrak g$.

Given a comultiplication $\delta \colon \mathfrak{g}\mapsto \mathfrak{g}\otimes \mathfrak{g}$ on a Lie algebra $\mathfrak{g}$, consider a vector space $D(\mathfrak{g})=\mathfrak{g}\oplus \mathfrak{g}^*$. 
As it was noted in \cite{Drinf}, a pair $(\mathfrak{g},\delta)$ is a Lie bialgebra if and only if there is a (unique) structure of a Lie algebra on $D(\mathfrak{g})$ satisfying

\noindent--- $\mathfrak{g}$ and $\mathfrak{g}^*$ are subalgebras in $D(\mathfrak{g})$, 

\noindent--- a bilinear form $Q \colon D(\mathfrak{g})\times D(\mathfrak{g})\mapsto F$ defined as $Q(x+f,y+g)=f(y)+g(x)$ for all $x,y\in \mathfrak{g}$, $f,g\in \mathfrak{g}$ is invariant.

In this case, the product in $D(\mathfrak{g})$ is given by the formula
$$
[x+f,y+g]_D=([x,y]+\langle f,y_{(2)}\rangle y_{(1)}+\langle g, x_{(1)}\rangle x_{(2)})+([f,g]+x\bullet g+f\bullet y),
$$
where  $\langle x\bullet g, y\rangle=\langle g\bullet y, x\rangle =g([y,x])$ for all $x,y\in \mathfrak{g}$, $f,g\in \mathfrak{g}^*$. The algebra $D(\mathfrak{g})$ is called {\it the classical double \cite{Stolin}} (or the Drinfeld double \cite{Zhel}) of the bialgebra $(\mathfrak{g},\delta)$.

Let $\mathfrak g$ be a Lie algebra and
$r=\sum\limits_ia_i\otimes b_i\in \mathfrak g\otimes \mathfrak g$. Define a comultiplication $\Delta_r$ on $\mathfrak g$ by
$$
\delta_r(a)=[r,a]
$$
for all $a\in \mathfrak g$. It is easy to see that $\delta_r$ is a 1-cocycle.
In \cite{BD} it was proved that $(\mathfrak g,\Delta_r)$ is a Lie coalgebra if
and only if the following two conditions hold:
\begin{itemize}
    \item The symmetric part $r+\tau(r)$ of the element $r$ is $\mathfrak g$-invariant.
    \item The element $C_{\mathfrak g}(r)=[r_{12},r_{13}]-[r_{23},r_{12}]+[r_{13},r_{23}]$
is $\mathfrak g$-invariant.
\end{itemize} 
Here $[r_{12},r_{13}]=\sum\limits_{ij}[a_i,a_j]\otimes b_i\otimes b_j$,
$[r_{23},r_{12}]=\sum\limits_{ij}a_i\otimes[a_j,b_i]\otimes b_j$,
 and $[r_{13},r_{23}]=\sum\limits_{ij} a_i\otimes a_j\otimes [b_i,b_j]$.

The equation 
\begin{equation}\label{CYBE}
 C_{\mathfrak g}(r)=[r_{12},r_{13}]-[r_{23},r_{12}]+[r_{13},r_{23}]=0
\end{equation}

is called the classical Yang-Baxter equation on $\mathfrak g$. A solution $r\in \mathfrak g\otimes \mathfrak g$ is called a classical $r$-matrix.

Let $\mathfrak{g}$ be a finite-dimensional Lie algebra over a field $F$ and $\bar{F}$ be an arbitrary field extension of $F$ and $\delta$ be a comultiplication. Consider $\mathfrak{g}_{\overline{F}}=\mathfrak{g}\otimes_F \overline{F}$ as a Lie algebra over $\overline{F}$ and define

 $$\delta_{\overline{F}}=\delta\otimes id_{\overline{F}} :\mathfrak{g}_{\overline{F}} \mapsto (\mathfrak{g} \otimes_F \mathfrak{g})\otimes _F \overline{F}\cong (\mathfrak{g} \otimes_F \overline{F})\otimes_{\overline{F}} (\mathfrak{g}\otimes _F \overline{F})=\mathfrak{g}_{\overline{F}}\otimes _{\overline{F}} \mathfrak{g}_{\overline{F}}.$$ 

It is straightforward to check that the classical double $D(\mathfrak{g}_{\overline{F}})$ of the bialgebra $(\mathfrak{g}_{\overline{F}},\delta_{\overline{F}})$ is isomorphic to $D(\mathfrak{g})\otimes_F \overline{F}$. Thus, if $(\mathfrak{g},\delta)$ is a Lie bialgebra, then $(\mathfrak{g}_{\overline{F}},\delta_{\overline{F}})$ is also a Lie bialgebra.

{\bf Definition.} Let $(\mathfrak{g},\delta)$ be a Lie bialgebra. Then

(1) The pair $(\mathfrak{g},\delta)$ is called a {\it coboundary} Lie bialgebra, if $\delta=\delta_r$ for some $r\in \mathfrak{g} \wedge \mathfrak{g}$.

(2) The pair $(\mathfrak g,\delta)$ is called a {\it triangular} Lie bialgebra, if $\delta=\delta_r$, for some $r\in \mathfrak g\wedge \mathfrak g$, where $r$ is a solution of CYBE on $\mathfrak{g}$.

(3)  The pair $(\mathfrak{g},\delta)$ is called a {\it quasitriangular} Lie bialgebra, if $\delta=\delta_r$ for some $r\in \mathfrak{g}\otimes \mathfrak{g}$, where $r$ is a~solution of CYBE on $\mathfrak{g}$.

(4) The pair $(\mathfrak{g},\delta)$ is called a {\it factorizable} Lie bialgebra, if $\delta=\delta_r$ for some $r\in \mathfrak g\otimes\mathfrak g$, where $r$ is a solution of CYBE on $\mathfrak g$ and $r+\tau(r)$ is non-degenerate (as a bilinear form on $\mathfrak{g}$) (\cite{STSRESH}).

(5) A Lie bialgebra $(\mathfrak{g},\delta)$ over a field $F$ is called {\it almost-factorizable}, if for some field extension $\bar{F}$, the Lie bialgebra $(\mathfrak{g}\otimes_F \overline{F},\delta_{\overline{F}})$ is factorizable (\cite{AndJan}).

It is known that if $\mathfrak{g}$ is a simple complex finite-dimensional Lie algebra, then any structure of a Lie bialgebra is either triangular or factorizable, that is, is induced by a special type Rota-Baxter operator of a weight $\lambda\in \mathbb C$.  If $\mathfrak{g}$ is an almost simple real finite-dimensional Lie algebra (that is, the complexification of $\mathfrak{g}$ is a complex simple Lie algebra), then in \cite{AndJan} it was noted that a Lie bialgebra structure on $\mathfrak{g}$ is either triangular, factorizable, or almost-factorizable (a similar result for reductive Lie algebras, that are extensions of finite-dimensional almost-simple real Lie algebras, follows from  \cite{FarJan}).

\section{Rota-Baxter operators of a non-scalar weights.}

Let $A$ be an arbitrary algebra over a field $F$. Given $x\in A$, denote by $L_x$ and $R_x$ operators of left and right multiplications:
$$
L_x(y)=xy,\ \ R_x(y)=yx
$$
for all $y\in A$. Clearly, $L_x,R_x\in \mathrm{End}(A)$. Denote by $\mathfrak{U}$ the subalgebra in $\mathrm{End}(A)$, generated by operators $L_x$ and $R_x$ for all $x\in A$: $\mathfrak U=alg\langle L_x,R_x|x\in A\rangle$. 

{\bf Definition \cite{Jac}.} A linear map $\mu:A\rightarrow A$ is called a centrum map in (a centralizer) of an algebra $A$ if, for any $x\in A$, $\mu\circ L_x=L_x\circ \mu$ and $\mu\circ R_x=R_x\circ \mu$. In other words,  for all $x,y\in A$
$$
\mu(xy)=\mu(x)y=x\mu(y).
$$

We will denote by ${\rm Cent}(A)$ the set of all centrum maps of the algebra $A$. In other words, $\mathrm{Cent}(A)$ is the centralizer of the subalgebra $\mathfrak U$ in $\mathrm{End}(A)$.  Clearly, ${\rm Cent}(A)$ is a subalgebra in ${\rm End}(A)$.

Note that if $A^2=A$, then the algebra ${\rm Cent}(A)$ is commutative. Indeed, for any $\mu,\gamma\in {\rm Cent}(A)$ and $x,y\in A$, we have
$$
(\varphi\circ\mu)(xy)=\varphi(x\mu(y))=\varphi(x)\mu(y)=(\mu\circ\varphi)(xy).
$$

{\bf Example 1.} Let $A$ be an arbitrary algebra over a field $F$ and $\lambda \in F$. Define a map $\varphi_{\lambda}$ as
\begin{equation}\label{Ex1}
\varphi_{\lambda}(x)=\lambda x
\end{equation}
for all $x\in A$. Obviously, $\varphi_{\lambda}\in {\rm Cent(A)}$ . 

{\bf Example 2.} Let $\mathfrak g=sl_2(\mathbb C)$ be the complex 3-dimensional special Lie algebra with the standard basis $x,h,y$ and the multiplication
$$
[x,y]=h,\ [h,x]=2x,\ [h,y]=-2y.
$$
Consider $\mathfrak g$ as a Lie algebra over $\mathbb R$. Then $\mathfrak g$ is a simple 6-dimensional Lie algebra and
$$
\mathfrak g=sl_2(\mathbb R)\oplus i sl_2(\mathbb R).
$$

Consider a map $\varphi_i$ defined as
$$
\varphi_i(a+ib)=i(a+ib)=ia-b$$
for all $a,b\in sl_2(\mathbb R)$. Then $\varphi_i\in {\rm Cent(\mathfrak g)}$ (see \cite{Goto} for details).  Moreover, a map $\psi$, defined as
$$
\psi(a+ib)=a-ib
$$
for all $a,b\in sl_2(\mathbb R)$, is an automorphism of $\mathfrak g$ (as an algebra over $\mathbb R$). Simple calculations show that 
$$\psi^{-1}\circ\varphi_i\circ \psi=-\varphi_i.$$

{\bf Example 3.} Suppose that an algebra $A$ splits into a direct sum of ideals $A_i(i=1,2,3)$: $A =A_1\oplus A_2\oplus A_3$. Consider a map $\mu$ defined as
$$
\mu(x_1+x_2+x_3)=x_2+2x_3,
$$
where $x_i\in A_i(i=1,2,3)$. Then $\mu$ is a centrum map with $ker(\mu)=A_1$.

{\bf Definition.} Let $A$ be an arbitrary algebra and $\mu\in \mathrm{Cent}(A)$. A linear map $R:A \rightarrow A$ is called a Rota-Baxter operator of weight $\mu$ if for all $x,y\in A:$
\begin{equation}\label{RB}
    R(x)R(y)=R(R(x)y+xR(y)+\mu (xy)).
\end{equation}

If $R$ is a Rota-Baxter operator on $A$, then we will call the pair $(A,R)$ a Rota-Baxter algebra.

{\bf Remark 1.} If the weight $\mu$ of a Rota-Baxter operator $R$ is equal to $\varphi_{\lambda}$ (defined by \eqref{Ex1} for some $\lambda\in F$), we will say that $R$ is a Rota-Baxter operator of weight $\lambda$ or that $R$ is a Rota-Baxter operator of a scalar weight.

{\bf Example 4.} Let $A$ be the algebra of continuous functions on $\mathbb R$. For any $f(x)\in A$, define
$$
R(f(x))=\int\limits_0^x f(t)dt.
$$

Then $R$ is a Rota-Baxter operator of weight 0 (due to the integration by parts formula) \cite{Baxter}.

{\bf Example 5.} Let $d$ be an inevitable derivation of an arbitrary  algebra $A$. Then $d^{-1}$ is a Rota-Baxter operator of weight 0.

The following statements generalize well-known results for Rota-Baxter operators of scalar weights to the general case.

{\bf Statement 1.} Suppose an algebra $A$ splits into a direct sum of two subalgebras:$A=A_1\oplus A_2$.  For any $\mu\in \mathrm{Cent}(A)$ such that $A_1^2$ and $A_2^2$ are $\mu$-invariant, a map $R$ defined as
$$
R(x+a)=-\mu (x)
$$
for all $x\in A_1$, $a\in A_2$, is a Rota-Baxter operator of weight $\mu$.

{\bf Proof.} Let $x,a\in A_1$, $y,b\in A_2$. Then
$$
R(x+a)R(y+b)=\mu(x)\mu(y)=\mu^2(xy).
$$
On the other hand,
\begin{multline*}
R(R(x+a)(y+b)+(x+a)R(y+b)+\mu((x+a)(y+b)))\\
=R(-\mu(x)(y+b)-(x+a)\mu(y)+\mu(x)(y+b)+a\mu(y+b))\\
=R(-x\mu(y)+a\mu(b))=R(-\mu(xy)+\mu(ab))=\mu^2(xy).
\end{multline*}

{\bf Statement 2.} Let $R$ be a Rota-Baxter operator of weight $\mu$.\\
(1) If $\gamma\in {\rm Cent}(A)$, then $R'=R\circ \gamma$ is a Rota-Baxter operator of weight $\mu\circ\gamma$. \\
(2) If $\varphi$ is an automorphism or an anti-automorphism of $A$, then $\varphi\circ R\circ \varphi^{-1}$ is a Rota-Baxter operator of weight $\varphi\circ \mu\circ\varphi^{-1}$.\\
(3) A map $R_1=-\mu-R$ is also a Rota-Baxter operator of the same weight $\mu$.\\
(4) If a linear map $Q:A\rightarrow A$ satisfies
$$
Q(x)Q(y)-Q(Q(x)y+xQ(y))=-\mu^2([x,y])
$$
for some $\mu\in\mathrm{Cent(A)}$ and all $x,y\in A$, then a map $R=\frac{1}{2}(Q-\mu)$ is a Rota-Baxter operator of weight $\mu$.

{\bf Proof.}
(1) For any $x,y\in A$, we have
\begin{multline*}
R'(x)R'(y)=R(\gamma(x))R(\gamma(y))=R(R(\gamma(x))\gamma(y)+\gamma(x)R(\gamma(y))+\mu (\gamma(x)\gamma(y)))\\
=R(\gamma(R'(x)y+xR'(y)+\gamma(x)\mu(\gamma(y)))=R'(R'(x)y+yR'(x)+(\mu\circ\gamma)(xy)).
\end{multline*}

(2) Similar to (1).

(3) For any $x,y\in\mathfrak g$, we have
\begin{multline*}
R_1(x)R_1(y)=(\mu+R)(x)(\mu+R)(y)\\
=\mu^2(xy)+\mu(R(x)y+xR(y))+R(x)R(y).
\end{multline*}

On the other hand,
\begin{multline*}
R_1(R_1(x)y+xR_1(y)+\mu(xy))\\
=R_1(-\mu(x)y-R(x)y-x\mu(y)-xR(y)+\mu(xy))\\
=(\mu+R)(R(x)y+xR(y)+\mu(xy))=\mu(R(x)y+xR(y)+\mu(xy))+R(x)R(y).
    \end{multline*}

    (4) similar to (3). The statement is proved.

{\bf Remark 2.} Results from Statement 2 were obtained in \cite{Baxter} in the case when $\mu=\varphi_{\lambda}$, $\lambda\in F$.

If $R:\mathfrak g\rightarrow \mathfrak g$ is a Rota-Baxter operator of a scalar weight $\lambda$, a classical result says that a pair $(\mathfrak g, [,]_R)$, where 
$$
[x,y]_R=[R(x),y]+[x,R(y)]+\lambda [x,y],\ x,y\in \mathfrak g,
$$
is also a Lie algebra \cite{STS}. A similar technique allows one to obtain a similar result in the case when the weight of $R$ is an arbitrary $\mu\in\mathrm{Cent}(\mathfrak g)$ (and arbitrary varieties of algebras). Nevertheless, we will use the following construction that we find interesting and that we will use in the future. This construction (with a little difference) originally appears  in \cite{GonGub} for Rota-Baxter operators of scalar weights and generalizes the construction of the classical double for Lie bialgebras. 

Let $A$ be an algebra over an arbitrary field $F$, $R\colon A \mapsto A$ be a linear map, and $\mu\in \mathrm{Cent}(A)$. Consider a~direct sum of vector spaces
$D_{R,\mu}(A)=A\oplus \bar{A}$, where $\bar{A}$ is an isomorphic copy of $A$. Define a product on $D_{R,\mu}(A)$ as follows:
\begin{multline*}
(a+\bar{b})*(x+\bar{y})= ax+R(ay)-aR(y)+R(bx)-R(b)x
\\+\overline{ay+bx-R(b)y-bR(y)-\mu(by)},
\end{multline*}
$a,b,x,y\in A$.

Obviously, $A$ is a subalgebra in $D_{R,\mu}(A)$.

For any $x\in A$, define
\begin{equation} \label{i-map}
i(x)=\bar{x}+(\mu(x)+R(x))\in D_{R,\mu}(A). 
\end{equation}
Define
\begin{equation}\label{I(G)}
I(A)=\{i(x)\mid x\in A\}
\end{equation}
Clearly, $I(A)$ is a subspace of $D_{R,\mu}(A)$ isomorphic to $A$ as a vector space, and $I(A)\cap A=0$. Moreover, 
$$
D_{R,\mu}(A)=A\oplus I(A)
$$
as a vector space. Also, for any $x,y\in A$, we have
\begin{multline}\label{s31}
i(x)*y
 = (\bar{x}+\mu(x)+R(x))*y
 \\= \overline{xy}+R(xy)-R(x)y+\mu(x)y+R(x)y
 = i(xy),
\end{multline}
that is, $A*I(A)\subset I(A)$. Similarly, $y*i(x)=i(yx)$ and $I(A)*A\in I(A)$.

{\bf Statement 3.} $I(A)$ is an ideal in $D_{R,\mu}(A)$ if and only if $R$ is a Rota-Baxter operator of weight $\mu$.

{\bf Proof.} 
In light of \eqref{s31}, we need to check when $i(x)*\overline{y}\in I(A)$ and $\overline{y}*i(x)\in I(A)$. Consider $x,y\in A$, then

\vspace{-0.9cm}
\begin{multline*}
i(x)*\bar{y}
 = (\bar{x}+\mu(x)+R(x))*\bar{y} \\
 = \overline{-\mu(xy) -xR(y)-R(x)y+\mu(x)y+R(x)y} \\
   +R(\mu(x)y)-\mu(x)R(y)+R(R(x)y)-R(x)R(y) \\
 = \overline{-xR(y)}-\mu(xR(y))-R(xR(y))+R(\mu(xy)+R(x)y+xR(y))-R(x)R(y)\\
 =i(-xR(y))+R(\mu(xy)+R(x)y+xR(y))-R(x)R(y).
\end{multline*}
Since $I(A)\cap A=0$, we obtain that $i(x)*\overline{y}$ belongs to $I(A)$ for any $x,y\in A$ if and only if $R$ is a Rota-Baxter operator of weight $\mu$. Similarly, $\bar{y}*i(x)\in I(A)$ if and only if $R$ is a Rota-Baxter operator of weight $\mu$. The statement is proved.

Note that if $R$ is a Rota-Baxter operator of weight $\mu$, we have for all $x,y\in A$
\begin{multline}\label{s32}
i(x)*i(y)
 = i(x)*(\bar{y}+\mu(y)+R(y))
 \\ = -i(xR(y))+i(\mu(xy+xR(y))
 = i( \mu(xy)).
\end{multline}
If $\mu$ is invertible, then a map
$i_{\mu}\colon x\mapsto i(\mu^{-1}(x))$ is an isomorphism of algebras $A$ and $I(A)$. If $\mu=0$, then $I(A)^2=0$.

{\bf Corollary 1.} Let $A$  be an algebra over a field $F$ and $R:A\rightarrow A$ be a Rota-Baxter operator of weight $\mu\in \mathrm{Cent}(A)$. Then a non-associative polynomial $f(t_1,\ldots,t_n)$ in non-commuting variables $t_1,\ldots,t_n$ is an identity in $A$  if and only if $f(t_1,\ldots,t_n)$ is an identity in $(D_{R,\mu}(A),*)$. In particular, algebras $(D_{\mu,R(A)},*)$ and $A$ belong to the same variety of algebras.

{\bf Proof.} Recall that $D_{R,\mu}(A)=A\oplus I(A)$. Let $f(t_1\ldots t_n)$ be a non-associative polynomial in non-commuting variables $t_1,\ldots t_n$. Suppose that $f(t_1,\ldots,t_n)$ is an identity in $A$, that is, for all $a_1,\ldots,a_n\in A$: $f(a_1,\ldots,a_n)=0$. Consider homogeneous elements $h_i\in A\cup I(A)$, $i=1,\ldots,n$. If for all $i=1,\ldots n$ $h_i\in A$, then $f(h_1,\ldots,h_n)=0$ by the assumption. In general, we have that $f(h_1,\ldots,h_n)=i(\mu^k(f(x_1,\ldots,x_n))=0$, where $k$ is a number of elements $h_1,\ldots h_n$ from $I(A)$ and $x_s=h_s$ if $h_s\in A$ and $x_s=i^{-1}(h_s)$ if $h_s\in I(A)$. 

Since $A$ is a subalgebra in $D_{R,\mu}(A)$, $A$ satisfy any identity of the algebra $D_{R,\mu}(A)$. The corollary is proved.

{\bf Corollary 2.} Let $\mathcal{M}$ be a variety of algebras, and $A\in\mathcal{M}$. Suppose that $R$ is a Rota-Baxter operator of weight $\mu\in\mathrm{Cent}(A)$. Define a new product $\cdot_R$ by
$$
x\cdot_R y=R(x)y+xR(y)+\mu(xy),
$$
where $x,y\in A$. Then the algebra $(A,\cdot_R)$ lies in the variety generated by the algebra $A$. In particular, $(A,\cdot_R)$ is an algebra from $M$.

{\bf Corollary 3.} Let $A$ be an algebra over an arbitrary field $F$ and $R:A\rightarrow A$ be a Rota-Baxter operator of weight $\mu\in\mathrm{Cent}(A)$. Then\\
(1) If $I(A)$ is the ideal in $D_{R,\mu}(A)$ defined in \eqref{I(G)}, then $I(A)^2=0$ if and only if $R$ is a Rota-Baxter operator of weight 0.\\
(2) If $\mu$ is an invertible map, then $D_{R,\mu}(A)$ is isomorphic to the direct sum  of ideals $A\oplus A$.

{\bf Proof.} 
From \eqref{s32}, it follows that $I(A)$ is an abelian ideal if and only if $\mu(xy)=0$ for all $x,y\in A$. Thus, the summand $\mu(xy)$ may be omitted and $R$ is a Rota-Baxter operator of weight 0.

Consider the second statement. If $\mu$ is invertible, then a map $\phi:A\mapsto I(A)$ defined as 
$$
\psi(x)=i(\mu^{-1}(x))
$$
is an isomorphism of algebras. Consider $J(A)=\{j(x)=\mu(x)-i(x)|x\in A\}$. From the definition of $I(A)$, we get that $I(A)\cap A=0$. Therefore, $I(A)\cap J(A)=0$.  For any $x\in A$
$$
x=i(\mu^{-1}(x))+j(\mu^{-1}(x)).
$$
Thus, $D_R(A)=I(A)\oplus J(A)$. From \eqref{s31} we get that
$$
j(x)*y=(\mu(x)-i(x))*y=\mu(xy)-i(xy)=j(xy)\in J(A)
$$
for any $x,y\in A$. From \eqref{s31} and \eqref{s32}, we obtain
$$
j(x)*i(y)=(\mu(x)-i(x))*i(y)=i(\mu(xy))-i(\mu(xy))=0
$$
for any $x,y\in A$. Simirarely, $y*j(x)\in J(A)$ and $i(y)j(x)=0$. Therefore, $J(A)$ is an ideal in $D_{R,\mu}(A)$, and $I(A)*J(A)=J(A)*I(A)=0$. Finally, for any $x,y\in A$,
$$
j(x)*j(y)=(\mu(x)-i(x))*j(y)=j(\mu(xy)).
$$
Thus, algebras $I(A)$, $J(A)$ and $A$ are isomorphic. The corollary is proved.
We will also need the following

{\bf Proposition 1.} Let $A$ be a simple algebra, $R\in\mathrm{End}(A)$ and $\mu\in \mathrm{End}(A)$. If $I$ is a proper ideal in $D_{R,\mu}$, then $I$ is isomorphic to $A$ as a vector space, and $D_{R,\mu}(A)=A\oplus I$ (a semidirect sum).

{\bf Proof.} Let $I$ be a proper ideal in $D_{R,\mu}(A)$. Then the intersection $I_1=I\cap A$ is an ideal in $A$. Therefore, $I_1=A$ or $I_1=0$. 

If $I_1=A$, then for any $a,b\in A$, we have that 
$$
a*\overline{b}=\overline{ab}+x
$$
for some $x\in A$. Since $A^2=A$ and $A\subset I$, we get that $\overline{A}\subset I$. Thus, $I=D_{R,\mu}(A)$, a contradiction.

Therefore, $A\cap I=0$. Let $W=\{a\in A|\ \overline{a}+b\in I\ \text{for some}\ b\in A\}$. Consider $a\in W$ and  $x\in A$. Suppose that $\overline{a}+b\in I$, where $b\in A$. Then 
$$
(\overline{a}+b)*x=\overline{ax}+y
$$
for some $y\in A$. Thus, $W$ is an ideal in $A$. If $W=0$, then $I=0$. Therefore, $W=A$. If for some $a\in A$ there are $b_1,b_2\in A$ such that $\overline{a}+b_1\in I$ and $\overline{a}+b_2\in I$, then $b_1-b_2\in I\cap A=0$. Thus, there is a linear function $\gamma\in\mathrm{End}(A)$ such that $I=\{\overline{a}+\gamma(a)|a\in A\}$. A map $\varphi: A\rightarrow I$ acting as
$$
\varphi(a)=\overline{a}+\gamma(a)
$$
for $a\in A$, is an isomorphism of vector spaces. 

Since $\overline{A}\subset A\oplus I$, we have that $D_{R,\mu}(A)=A\oplus I$. The proposition is proved.

\section{Rota-Baxter operators and CYBE}

Rota-Baxter operators (of a scalar weight) naturally appear in the paper of M.A. Semenov-Tyan-Shanskii \cite{STS}. In this paper (titled "What a classical $r$-matrix really is", a "right" version of the notion of a classical $r$ matrix was suggested. Instead of considering the notion of Lie bialgebra, he introduced the notion of a double Lie algebra (which should not be confused with the notion of the double Lie algebra from \cite{Sole} and \cite{GonGub}). 

{\bf Definition (\cite{STS}).} Given a Lie algebra $\mathfrak{g}$, a linear map $R\in \textrm{End}(\mathfrak{g})$ is called a classical $R$-matrix if the space $\mathfrak{g}$ with a new product $[\cdot,\cdot]_R$, defined as
$$
[x,y]_R=[R(x),y]+[x,R(y)],
$$
where $x,y\in \mathfrak{g}$, is a Lie algebra.

{\bf Remark 3.} In \cite{STS}, it was noticed that the notion of the double Lie algebra coincides with the notion of the Lie bialgebra if and only if the map $R$ is skew-symmetric with respect to some non-degenerate bilinear quadratic form, see Section 3 for details.

Given a map $R:\mathfrak g\mapsto \mathfrak g$, define a bilinear map $\theta_R$ as
$$
\theta_R(x,y)=[R(x),R(y)]-R([R(x),y]+[x,R(y)]).
$$

{\bf Proposition 2 (\cite{STS}})
$R$ is a classical R-matrix if and only if 
\begin{equation}\label{e4}
[\theta_R(x,y),z]+[\theta_R(y,z),x]+[\theta_R(z,x),y]=0,\ \forall x,y,z\in \mathfrak{g}.
\end{equation} 

M.A Semenov-Tan-Shanskii suggested the following two conditions, each of which implies \eqref{e4}: 

i) $\theta_R(x,y)=0$.

In this case, $R$ is a Rota-Baxter operator of weight zero. In \cite{STS}, this equation is also called the {\it classical Yang--Baxter equation}.

ii) $\theta_R(x,y)=- \lambda^2 [x,y]$, $\lambda\neq 0$. 

This equation is called the {\it modified classical Yang--Baxter equation} (mCYBE).  It can be written as
\begin{equation}\label{e5}
[R(x),R(y)]=R([R(x),y]+[x,R(y)])-\lambda^2[x,y].
\end{equation}

 The parameter was chosen so the map $R=\lambda \textrm{id}$ is a solution of \eqref{e5}.

Conditions (i) and (ii) are easier to work with since they are bilinear rather than trilinear.
Note that the condition $\theta_R([x,y])=\mu([x,y])$ for some $\mu\in \mathrm{Cent}(\mathfrak g)$ also implies \eqref{e4}.

It turns out that a map $R\in \textrm{End}(\mathfrak{g})$ is a solution of mCYBE \eqref{e5} 
if and only if  a map
$$
B=\frac{1}{2}(R-\lambda \textrm{id}) 
$$
is a Rota--Baxter operator of weight $\lambda$ on $\mathfrak{g}$.

{\bf Definition} Let $\mathfrak{g}$ be a Lie algebra. A bilinear form
$\omega$ on $\mathfrak{g}$ is called invariant if for all $a,b,c\in \mathfrak{g}$:
$$
\omega([a,b],c)=\omega(a,[b,c]).
$$

{\bf Definition.}
Let $\mathfrak{g}$ be a Lie algebra, and $\omega$ is a bilinear invariant
non-degenerate form on $\mathfrak{g}$. Then the pair $(\mathfrak{g},\omega)$ is called a
quadratic Lie algebra.

{\bf Notation.} Let $(\mathfrak{g},\omega)$ be a quadratic Lie algebra and $R:\mathfrak{g}\rightarrow \mathfrak{g}$ be a linear map. Then by $R^*$ we will denote the adjoint to $R$ with respect to the form $\omega$ map defined as
$$
\omega(R(x),y)=\omega(x,R^*(y))
$$
for all $x,y\in \mathfrak{g}$.

{\bf Remark 4.} Let $(\mathfrak g,\omega)$ be a quadratic Lie algebra and $\mu\in \mathrm{Cent}(\mathfrak g)$. Then for all $x,y\in \mathfrak g$, $\mu^*([x,y])=\mu([x,y])$. Indeed, for all $x,y,z\in \mathfrak g$, we have
$$
\omega(\mu([x,y]),z)=\omega([z,\mu(y)],z)=\omega(z,[\mu(y),z])=\omega(z,[y,\mu(z)])=\omega([z,y],\mu(z)).
$$

In particular, if $\mathfrak g$ is perfect, then  $\mu^*=\mu$.

If $(\mathfrak{g},\omega)$ is a quadratic finite-dimensional Lie algebra, then there is a natural isomorphism $^* \colon\mathfrak{g}\mapsto \mathfrak{g}^*$ between the vector space $\mathfrak{g}$ and the dual space $\mathfrak{g}^*$ defined as $x^*(y)=\omega(x,y)$. This, in turn, gives a natural isomorphism between $\textrm{End}(\mathfrak{g})\cong \mathfrak{g}\otimes \mathfrak{g}^*$ and $\mathfrak{g}\otimes \mathfrak{g}$. This allows us to consider an element $r=\sum a_i\otimes b_i\in \mathfrak{g}\otimes \mathfrak{g}$ as a linear map $R_r$ acting as
\begin{equation}\label{TenzOp}
R_r(x)=\sum \omega (a_i,x) b_i.
\end{equation}

{\bf Remark 5.} Note that the adjoint map $R_r$ with respect to the form $\omega$ is defined by
$$
R^*_r(x)=\sum \omega (b_i, x) a_i
$$
for all $x\in \mathfrak{g}$.

{\bf Statement 4 (\cite{STS}).} Let $(\mathfrak g, \omega)$ be a finite-dimensional quadratic Lie algebra. An element $r\in \mathfrak g\otimes \mathfrak g$ is a skew-symmetric solution of the classical Yang-Baxter equation \eqref{CYBE} if and only if the corresponding map $R_r$ is a skew-symmetric Rota-Baxter operator of weight 0.

In other words, skew-symmetric solutions of the classical Yang-Baxter equation \eqref{CYBE} are in one-to-one correspondence with skew-symmetric Rota-Baxter operators in the case of quadratic Lie algebras.

Consider a structure of a factorizable Lie bialgebra $(\mathfrak g,\delta_r)$ on a Lie algebra $\mathfrak g$. Since $I=r+\tau(r)$ is non-degenerate, it defines an isomorphism of vector spaces $I:\mathfrak g^*\mapsto \mathfrak g$: $I(f)=(f\otimes id)(r)$. In \cite{STSRESH} it was noted that a form $\beta$ defined as
\begin{equation}\label{RBform}
\beta(x,y)=\langle I^{-1}(x), y\rangle
\end{equation}
for all $x,y\in \mathfrak g$ is a nondegenerate invariant form on $\mathfrak g$. Moreover,  in \cite{BD2} it was proved that the corresponding to $r$ and $\beta$ linear map $R_r$  is a Rota-Baxter operator of weight 1 satisfying
\begin{equation}\label{RBfact}
R_r+R_r^*+id=0.
\end{equation}
In other words, there is a one-to-one correspondence between structures of a factorizable Lie bialgebras on a finite-dimensional Lie algebra $\mathfrak g$ and triples $(\mathfrak g, \beta, R)$, where $\beta$ is a non-degenerate invariant form on $\mathfrak g$ and $R$ is a Rota-Baxter operator of weight 1 satisfying \eqref{RBfact}. 

Given a finite-dimensional Lie algebra $\mathfrak g$, there may be more than one non-degenerate invariant form on $\mathfrak g$ (even if $\mathfrak g$ is simple). We can start with a quadratic finite-dimensional coboundary Lie bialgebra $(\mathfrak g,\delta_r)$, $r\in \mathfrak g\otimes \mathfrak g$ and ask when $(\mathfrak g,\delta_r)$ is factorizable. We obtain the following version of the result from \cite{STSRESH}.

{\bf Theorem 1.} Let $\mathfrak g$ be a finite-dimensional Lie algebra over a field $F$ and $\omega$ be a non-degenerate bilinear invariant form on $\mathfrak g$. Let $r=\sum\limits_i a_i\otimes b_i\in \mathfrak g\otimes \mathfrak g$ and $R_r$ be the corresponding map defined as in \eqref{TenzOp}. Then\\
(1) $r$ is a solution of CYBE with $\mathfrak g$-invariant symmetric part $r+\tau(r)$ if and only if the map $R_r$ is a Rota-Baxter operator of a weight $\mu\in\mathrm{Cent}(\mathfrak g)$ satisfying
\begin{equation}\label{l11}
R_r+R_r^*+\mu=0,
\end{equation}
where $R^*$ is the adjoint to $R$ with respect to the form $\omega$ map.\\
(2) $r$ is a solution of CYBE with $\mathfrak g$-invariant and non-degenerate symmetric part $r+\tau(r)$ if and only if the map $R_r$ is a Rota-Baxter operator of a non-degenerate weight $\mu\in\mathrm{Cent}(\mathfrak g)$ satisfying \eqref{l11}.

{\bf Proof.}  For $x,y\in \mathfrak g$, consider a map $\psi_{x,y}:\mathfrak g\otimes \mathfrak g\otimes
\mathfrak g\rightarrow \mathfrak g$ defined as
$$
\psi_{x,y}(a\otimes b\otimes c)=\omega(x,a)\omega(y,b)c
$$
for all $a,b,c\in \mathfrak g$.

 Since $\omega$ is non-degenerate, $C_{\mathfrak g}(r)=0$ if and only if $\psi_{x,y}(C_{\mathfrak g}(r))=0$ for all $x,y\in \mathfrak g$.  Direct computations show
$$
\psi_{x,y}(C_{\mathfrak g}(r))=\sum\limits_{i,j}\omega([a_i,a_j],x)\omega(b_i,y)b_j-\omega(a_i,x)\omega([a_j,b_i],y)b_j+\omega(a_i,x)\omega(a_j,y)[b_i,b_j]=
$$
$$
= R_r([x,R^*_r(y)]-[R(x),y])+[R_r(x),R_r(y)].
$$
That is, $r$ is a solution of CYBE if and only if for all $x,y\in \mathfrak g$
\begin{equation}\label{t11}
    R_r([x,R^*_r(y)]-[R(x),y])+[R_r(x),R_r(y)]=0.
\end{equation}
Similarly, $r+\tau(r)$ is $\mathfrak g$-invariant if and only if for all $x,y\in \mathfrak g$
$$
[(R_r+R_r^*)(x),y]=(R_r+R_r^*)([x,y]),
$$
that is, $R_r+R_r^*$ satisfies \eqref{l11} for some $\mu\in \mathrm{Cent}(\mathfrak g)$. Note that \eqref{t11} and \eqref{l11} are equivalent to the fact that $R_r$ is a Rota-Baxter operator of weight $\mu$ satisfying \eqref{l11}. Finally, $r+\tau(r)$ is non-degenerate if and only if $\mu=-R_r-R^*_r$ is a non-degenerate map. The theorem is proved.

{\bf Statement 5.} Let  $\mathfrak g$ --- be a Lie algebra over a field $F$ with $Z(\mathfrak g)=0$, $\omega:\mathfrak g\times \mathfrak g\mapsto F$ be a bilinear invariant non---degenerate form on $\mathfrak g$. Suppose that  $\theta :\mathfrak g\times \mathfrak g\mapsto \mathfrak g$ is a bilinear map satisfying the following conditions
\begin{gather}
       \theta(x,y)=-\theta(y,x), \label{ssym}\\
    \omega(\theta(x,y),z)=\omega(x,\theta(y,z)) \label{assoc}\\
    [\theta(x,y),z]=\theta([x,z],y)+\theta(x,[y,z]), \label{nleyb} 
\end{gather}
for all $x,y,z\in \mathfrak g$. Then there is a linear map $\nu: [\mathfrak g,\mathfrak g]\rightarrow \mathfrak g$ such that $\theta(x,y)=\nu([x,y])$ for all $x,y\in \mathfrak g$. Moreover, for all $x,y,z\in \mathfrak g$ the map $\nu$ satisfies
\begin{equation}\label{centroid}
\nu([[x,y],z])=[\nu([x,y]),z].
\end{equation}

{\bf Proof.}

Let $x,y,z\in \mathfrak g$. From \eqref{ssym} and \eqref{nleyb} it follows that
\begin{multline}\label{addeq}
[\theta(x,y),z]=\theta([x,z],y)+\theta(x,[y,z])=[x,\theta(z,y)]-\theta(z,[x,y])+[y,\theta(x,z)]-\theta([y,x],z)\\
=[x,\theta(z,y)]+[\theta(z,x),y]-2\theta (z,[x,y]).
\end{multline}

Therefore, for all $x,y,z,t\in \mathfrak g$:
\begin{multline}
    \omega\left ([\theta(x,y),z],t\right )=\omega\left ([x,\theta(z,y)]+[\theta(z,x),y]-2\theta (z,[x,y]),t\right)\\
    =\omega(z,\theta(y,[t,x])+\theta(x,[y,t]))-2\omega(\theta(z,[x,y]),t)
    =\omega(z,[t,\theta(y,x)])-2\omega(\theta(z,[x,y]),t)\\
    =\omega([z,\theta(x,y)]-2\theta(z,[x,y]),t).
\end{multline}

Since the form $\omega$ is non-degenerate, we have
\begin{equation}\label{meq}
    [\theta(x,y),z]=\theta([x,y],z).
\end{equation}

Consider $x=\sum\limits_i [x_i,y_i]\in [\mathfrak g,\mathfrak g]$, where $x_i,y_i\in \mathfrak g$. Define  
$$
\nu(x)=\sum\limits_i \theta(x_i,y_i).
$$

First, we need to prove that $\nu$ is well-defined. Indeed, if $\sum [x_i,y_i]=\sum [z_i,q_i]$ for some $x,y,z,q\in \mathfrak g$, then we can use \eqref{meq} and note that for any $z\in \mathfrak g$ we have:
$$
\sum[\theta(x_i,y_i),z]=\theta\left (\sum[x_i,y_i],z\right)=\theta\left(\sum [z_i,q_i],z\right)=\sum [\theta (z_i,q_i),z].
$$

Since the center of the algebra $\mathfrak g$ is equal to 0, we get that  $\sum \theta(x_i,y_i)=\sum \theta(z_i,q_i)$. Thus, $\nu$ is well-defined. Obviously,  $\nu$ is linear.

Consider $x,y,z\in \mathfrak g$. Then 
$$
\nu([[x,y],z])=\theta([x,y],z)=[\theta(x,y),z]=[\nu([x,y]),z].
$$ 
The statement is proved.

{\bf Corollary 4.} Let $\mathfrak g$ be a perfect finite-dimensional Lie algebra over a field  $F$ with $Z(\mathfrak g)=0$. Suppose that  $\theta :\mathfrak g\times \mathfrak g\mapsto \mathfrak g$ is a bilinear map satisfying conditions \eqref{ssym}-\eqref{nleyb}. Then there is a map $\nu\in\mathrm{Cent}(\mathfrak g)$ such that for all $x,y\in \mathfrak g$: 
$$
\theta(x,y)=\nu([x,y])=[\nu(x),y]=[x,\nu(y)].
$$

\section{Coboundary Lie bialgebra structures and Rota-Baxter operators.}

Let $(\mathfrak g,\omega)$ be a quadratic Lie algebra over a field different from 2, $r=\sum\limits_i a_i\otimes b_i\in \mathfrak g\wedge\mathfrak g$. Define a comultiplication $\delta_r$ as
$$
\delta_r(x)=[r,x]=\sum\limits_i [a_i,x]\otimes b_i+a_i\otimes [b_i,x]
$$
for all $x\in \mathfrak g$.

Consider a map $R_r:\mathfrak g\rightarrow \mathfrak g$ defined as
$$
R_r(x)=\sum\limits_i\omega(a_i,x)b_i
$$
for all $x\in \mathfrak g$. Since $r$ is skew-symmetric, $R^*=-R$. Recall that given $x\in \mathfrak g$, by $x^*$ we will denote the  element from the dual space $\mathfrak g^*$ defined as $\langle x^*,y\rangle=\omega(x,y)$ for all $y\in \mathfrak g$.

There are several questions that one can ask: when a pair $(\mathfrak g,\delta_r)$ is a Lie bialgebra? If $(\mathfrak g,\delta_r)$ is a Lie bialgebra, is it a triangular, factorizable, or  almost factorizable Lie bialgebra? What can we say about the classical double of $(\mathfrak g, \delta_r)$? In this section, we will consider these questions from the point of view of the map $R_r$.

{\bf Proposition 3.} Let $(\mathfrak g,\omega)$ be an arbitrary finite-dimentional quadratic Lie algebra, and let $R\colon \mathfrak{g}\rightarrow \mathfrak{g}$ be a linear operator satisfying $R+R^*+\mu=0$ for some $\mu\in \mathrm{Cent}(\mathfrak g)$ and $r\in \mathfrak{g} \otimes \mathfrak{g}$  be the corresponding to the map $R$ tensor defined by \eqref{TenzOp}. Then the multiplication in the dual algebra of the coalgebra $(\mathfrak{g},\delta_r)$ is given by
$$
[a^*,b^*]=(-[R(a),b]-[a,R(b)]-\mu([a,b]))^*.
$$
for all $a,b\in \mathfrak g$. Moreover, the classical double $D(g)=\mathfrak{g}\oplus \mathfrak{g}^*$ of the pair $(\mathfrak{g},\delta_r)$ (with the product defined by $[\cdot,\cdot]_D$) is isomorphic to the algebra $D_{R,\mu}(\mathfrak g)=\mathfrak{g}\oplus \overline{\mathfrak{g}}$ with a product $[\cdot,\cdot]_R$ given by
\begin{multline*}
[x+\bar{y},z+\bar{t}]_R=[x,z]+R([x,t])-[x,R(t)]+R([y,z])-[R(y),z] \\
 +\overline{[x,t]+[y,z]-[R(y),t]-[y,R(t)]-\mu( [y,t])}
\end{multline*}
for all $x,y,z,t\in \mathfrak{g}$.

{\bf Proof.}
Let $r=\sum a_i\otimes b_i$. Then $\delta_r(t)=\sum [a_i,t]\otimes b_i+a_i\otimes [b_i,t]$ for all $t\in \mathfrak{g}$. Consider a map $\varphi:D(\mathfrak{g})\rightarrow D_{r,\mu}(\mathfrak{g})$ defined as
$$
\varphi(x+y^*)=x+\overline{y}
$$
for all $x,y\in \mathfrak{g}$. Obviously, $\varphi$ is an isomorphism of vector spaces, and $\varphi([x,y]_D)=[\varphi(x),\varphi(y)]_R$ for all $x,y\in \mathfrak g$.

For any $t\in \mathfrak{g}$ we have
\begin{multline*}
\langle [x^*,y^*]_D, t\rangle=\omega(x,t_{(1)})\omega(y,t_{(2)})=\sum \omega (x,[a_i,t])\omega(y,b_i)+\omega(x,a_i)\omega(y,[b_i,t])\\
= \sum \omega ([t,x],a_i)\omega(y,b_i)+\omega(x,a_i)\omega([t,y],b_i)
=\omega([t,x],R^*(y))+\omega([t,y],R(x))\\
= \omega(t,-[x,R(y)]-[x,\mu(y)]-[R(x),y])=\langle (-[R(x),y]-[x,R(y)]-\mu([x,y]))^*, t\rangle.
\end{multline*}
Hence, $[x^*,y^*]_D=(-[R(x),y]-[x,R(y)]-\mu([x,y]))^*$ and $\varphi([x^*,y^*]_D)=[\varphi(x^*),\varphi(y^*)]_R$.

Consider $\varphi([x,y^*]_D)$ for some $x,y\in\mathfrak g$. For any $t\in \mathfrak{g}$, using the associativity of the form $Q$, we have the following equations:
\begin{multline}\label{p11}
Q([x,y^*]_D,t^*)=Q(x,[y^*,t^*]_D)=-\omega(x,[R(y),t]+[y,R(t)]+\mu([y,t]))\\
=\omega(-[x,R(y)],t)+\omega(x,[y,R^*(t)])=\omega(-[x,R(t)]+R([x,y]),t) 
\end{multline}
\begin{equation}\label{p12}
  Q([x,y^*]_D,t)=Q(y^*,[t,x])=\omega(y,[t,x])=\omega([x,y],t)=Q([x,y]^*,t).
\end{equation}
From \eqref{p11} and \eqref{p12} it follows that 
$$
[x,y^*]_D=R([x,y])-[x,R(y)]+[x,y]^*.
$$
 Therefore, $\varphi ([x,y^*]_D)=[\varphi(x),\varphi(y^*)]_R$, and $\varphi$ is an isomorphism of algebras. The proposition is proved.

From Theorem 1, Proposition 3, and Corollary 3,  it follows

{\bf Corollary 5.} Let $(\mathfrak g,\omega)$ be a quadratic finite-dimensional Lie algebra and $\delta:\mathfrak g\rightarrow \mathfrak g\otimes \mathfrak g$ be a comultiplication. Then\\
(1) If $(\mathfrak g,\delta)$ is a triangular  Lie bialgebra, then the classical double is a semidirect sum $D(\mathfrak g)=\mathfrak g\oplus D$, where $D$ is an abelian ideal in $D(\mathfrak g)$.\\
(2) If $(\mathfrak g,\delta)$ is a factorizable Lie bialgebra, then the classical double $D(\mathfrak g)$ is isomorphic to the direct sum of two ideals $\mathfrak g\oplus \mathfrak g$.

{\bf Definition.} Given a Lie algebra $\mathfrak g$ and a linear map $R:\mathfrak g \rightarrow \mathfrak g$, define a linear map $\theta: \mathfrak g\times \mathfrak g\rightarrow \mathfrak g$ as
\begin{equation}\label{p41}
\theta(x,y)=[R(x),R(y)]-R_r([R(x),y]+[x,R(y)]
\end{equation}
for all $x,y\in \mathfrak g$.

{\bf Proposition 4.} Let $(\mathfrak g,\omega)$ be a quadratic Lie algebra, $r\in \mathfrak g\wedge \mathfrak g$, $R_r$ be the corresponding map, and $\theta: \mathfrak g\times \mathfrak g\rightarrow \mathfrak g$ be a bilinear map defined by \eqref{p41}. Then for all $x,y,z\in \mathfrak g$ 
\begin{gather}
\label{p31}\theta(x,y)=-\theta(y,x),\\
\label{p32}\omega(\theta(x,y),z)=\omega (x,\theta(y,z)).
\end{gather}
{\bf Proof.} The equality \eqref{p31} is obvious. Let $x,y,z\in \mathfrak g$. Then
\begin{multline*}
    \omega(\theta(x,y),z)=\omega([R_r(x),R_r(y)]-R_r([R_r(x),y]+[x,R_r(y)]),z)\\
    =\omega(R(x),[R(y),z]-[y,R^*_r(z)])-\omega(x,[R_r(y),R^*(z)])\\
    =\omega(x,[R_r(y),R_r(z)]-R_r([R_r(y),z]+[y,R_r(z)])=\omega(x,\theta(y,z)).
\end{multline*}

{\bf Proposition 5.} Let $(\mathfrak g,\omega)$ be a quadratic Lie algebra, $r\in \mathfrak g\wedge \mathfrak g$, $R_r$ be the corresponding map, and $\theta: \mathfrak g\times \mathfrak g\rightarrow \mathfrak g$ be a bilinear map defined by \eqref{p41}. Then the pair $(\mathfrak g, \delta_r)$ is a Lie bialgebra if and only if the following condition holds
\begin{equation}\label{p43}
 [\theta(x,y),z]=\theta([x,z],y)+\theta(x,[y,z])
\end{equation}
for all $x,y,z\in \mathfrak g.$

{\bf Proof.} Since $\delta_r$ satisfies \eqref{bder} for any $r\in \mathfrak g\otimes \mathfrak g$, we need to check when the pair $(\mathfrak g, \delta_r)$ is a Lie coalgebra. From Proposition 3, we know that $(\mathfrak g,\delta_r)$ is a Lie coalgebra if and only if the pair $(\mathfrak g,[\cdot,\cdot])$, where the product is defined as
$$
[x,y]_{R_r}=[R_r(x),y]+[x,R_r(y)]
$$
for all $x,y\in \mathfrak g$, is a Lie algebra. From Proposition 2, it follows that  $(\mathfrak g,[\cdot,\cdot]_{R_r})$ is a Lie algebra if and only if the bilinear map $\theta$ defined in \eqref{p41} satisfies
\begin{equation}\label{p42}
[\theta(x,y),z]+[\theta(y,z),x]+[\theta(z,x),y]=0
\end{equation}
for all $x,y,z\in \mathfrak g$. Consider arbitrary $x,y,z,t\in \mathfrak g$. Then 
\begin{multline*}
    \omega([\theta(x,y),z]+[\theta(y,z),x]+[\theta(z,x),y],t)    =\omega(x, \theta(y,[z,t])+[t,\theta(y,z)]+\theta([y,t],z).
\end{multline*}
Therefore, conditions \eqref{p43} and \eqref{p42}
 are equivalent. The proposition is proved. 

 {\bf Theorem 2.} Let $(\mathfrak g,\omega)$ be a quadratic finite-dimensional Lie algebra, where $\mathfrak g$ is a perfect Lie algebra with $Z(\mathfrak g)=0$, $r\in \mathfrak g\wedge \mathfrak g$ be a skew-symmetric tensor, $R_r$ be the corresponding map, and $\theta: \mathfrak g\times \mathfrak g\rightarrow \mathfrak g$ be a bilinear map defined by \eqref{p41} Then $(\mathfrak g, \delta_r)$ is a Lie bialgebra if and only if $\theta(x,y)=\nu([x,y])$ for some $\nu\in \mathrm{Cent}(\mathfrak g)$.

{\bf Proof.} If $\theta(x,y)=\nu([x,y])$ for some $\nu\in\mathrm{Cent}(\mathfrak g)$, then $(\mathfrak g,\delta_r)$ is a Lie bialgebra by Proposition 4.

Suppose that $(\mathfrak g,\delta_r)$ is a Lie bialgebra. By Proposition 4 and Proposition 5, the map $\theta$ satisfies equations \eqref{ssym}-\eqref{nleyb}. By Statement 5, there is a map $\nu\in\mathrm{Cent}(\mathfrak g)$ such that $\theta(x,y)=\nu([x,y])$ for all $x,y\in \mathfrak g$. The theorem is proved.

{\bf Theorem 3.}  Let $(\mathfrak g,\omega)$ be a quadratic finite-dimensional Lie algebra, $r\in \mathfrak g\wedge \mathfrak g$ be a skew-symmetric tensor, $R_r$ be the corresponding map, and $\theta: \mathfrak g\times \mathfrak g\rightarrow \mathfrak g$ be a bilinear map defined by \eqref{p41} Then \\
(i) The pair $(\mathfrak g, \delta_r)$ is a triangular Lie bialgebra if and only if $\theta(x,y)=0$ for all $x,y\in \mathfrak g$.\\
(ii) The pair $(\mathfrak g,\delta_r)$ is a quasitriangular Lie bialgebra if and only if $\theta(x,y)=-\mu^2([x,y])$ for some $\mu\in \mathrm{Cent}(\mathfrak g)$ and all $x,y\in\mathfrak g$.\\
(iii) The pair $(\mathfrak g,\delta_r)$ is a factorizable Lie bialgebra if and only if $\theta(x,y)=-\mu^2([x,y])$ for some non-degenerate $\mu\in \mathrm{Cent}(\mathfrak g)$ and all $x,y\in\mathfrak g$.

{\bf Proof.} (i) This is equivalent to the fact that skew-symmetric solutions of the classical Yang-Baxter equations are in one-to-one correspondence with skew-symmetric Rota-Baxter operators of weight 0.

(ii) Suppose that $\theta(x,y)=-\mu^2([x,y])$ for some $\mu\in \mathrm{Cent}(\mathfrak g)$. 

First we note that we may assume that $\mu^*=\mu$. Indeed, we can consider $\nu=\frac{1}{2}(\mu+\mu^*)$. Then $\nu^2=\frac{1}{4}(\mu^2+\mu\mu^*+\mu^*\mu+(\mu^*)^2$). By Remark 4, $\mu([x,y])=\mu^*([x,y])$ for all $x,y\in \mathfrak g$. Thus,  $\mu^2([x,y])=\nu^2([x,y])$ for all $x,y\in \mathfrak g$ and we can take $\nu$ instead of $\mu$.

Suppose $\mu^*=\mu$. Consider a map $B_{r}$ defined as
$$
B_{r}(x)=\frac{1}{2}(R_{r}(x)-\mu(x))
$$
for all $x\in \mathfrak g$. By Statement 2, $B_r$ is a Rota-Baxter operator of weight $\mu$. Since $R_r+R_r^*=0$ and $\mu^*=\mu$, $B_r$ satisfies $B_r+B_r^*+\mu=0$. Let $r_1\in \mathfrak g\otimes \mathfrak g$ be the corresponding to the map $2B_r$ tensor \eqref{TenzOp}. By Theorem 1, $(\mathfrak g,\delta_{r_1})$ is a quasitriangular Lie bialgebra. 

 By Proposition 3, the multiplication in the dual algebra of the bialgebra $(\mathfrak g, \delta_r)$ is given by
$$
[a^*,b^*]=(-[R(a),b]-[a,R(b)])^*=(-[2B_r(a),b]-[a,2B_r(b)]-2\mu[a,b])^*
$$
for all $a,b\in \mathfrak g$. That is, the multiplication in the dual algebra of the coalgebra $(\mathfrak g,\delta_r)$ coincides with the product in the dual algebra of the coalgebra $(\mathfrak g, \delta_{r_1})$. Therefore, $\delta_r=\delta_{r_1}$ and $(\mathfrak g, \delta_r)$ is a quasitriangular Lie bialgebra. 

Suppose that $(\mathfrak g, \delta_r)$ is a quasitriangular Lie bialgebra. That is, $\delta_r=\delta_t$, where $t\in \mathfrak g\otimes \mathfrak g$ is a solution of the CYBE and $t+\tau(t)$ is non-degenerate. By Theorem 1, the corresponding to $t$ map $R_t$ is a Rota-Baxter operator of weight $\eta\in\mathrm{Cent}(\mathfrak g)$ and for all $x\in \mathfrak g$
$$
R_t(x)+R^*_t(x)+\eta(x)=0.
$$

Let $q=t-r=\sum\limits_i x_i\otimes y_i\in \mathfrak g\otimes \mathfrak g$. Since for any $x\in \mathfrak g$ $\delta_r(x)=\delta_t(x)$, we have that $[q,x]=0$ for any $x\in \mathfrak g$. Let $\mu:\mathfrak g\rightarrow \mathfrak g$ be the corresponding to $q$ map defined by \eqref{TenzOp}. We have for any $y\in \mathfrak g $
$$
0=(y^*\otimes id)([q,x])=\sum\limits \omega(y,[x_i,x])y_i+\omega(y,x_i)[y_i,x]=\mu([x,y])+[\mu(y),x].
$$
That is, $\mu\in\mathrm{Cent}(\mathfrak g)$. Therefore,
$$
R_t+R_t^*=R_r+R_r^*+\mu+\mu^*=\mu+\mu^*,
$$
that is, $\eta=-\mu-\mu^*$. Note that by Remark 4, $\nu([x,y])=2\mu([x,y])$ for all $x,y\in \mathfrak g$. 

Let $x,y\in \mathfrak g$. Since $R_t=R_r+\mu$ is a Rota-Baxter operator of weight $\eta$, for all $x,y\in \mathfrak g$ we have
\begin{multline*}
0=[R_t(x),R_t(y)]-R_t([R_t(x),y]+[x,R_t(y)]+\eta([x,y]))\\
=[(R_r(x)+\mu(x),R_r(y)+\mu(y)]-R_r([R_r(x)+\mu(x),y]+[x,R_r(y)+\mu(y)]+\eta([x,y]))\\
-\mu([R(x)+\mu(x),y]+[x,R_r(y)+\mu(y)]+\eta([xy])\\
=[R_r(x),R_r(y)]-R_r([R_r(x),y]+[x,R_r(y)])+\mu^2([x,y]).
\end{multline*}
Therefore, $\theta(x,y)=-\mu^2([x,y])$. The theorem is proved.

\section{Simple finite-dimensional Lie bialgebras.}
 
All the results in this section (except the last theorem) are not new in the case $F=\mathbb R$. Nevertheless, it may be convenient to have a unique approach to the topic for an arbitrary field of characteristic zero.

{\bf Definition.} An (arbitrary) algebra $A$ over a field $F$ is called  absolutely simple (normally simple in \cite{Jac}) if any field extension $F\subset \overline{F}$, the algebra $A\otimes _F\overline F$ is simple.

We will need the following classical characterization of absolutely simple algebras. 

{\bf Theorem 4\cite{Jac}.} An arbitrary finite-dimensional algebra $A$ is absolutely simple if and only if it is simple and $\mathrm{Cent}(\mathfrak g)=\{\lambda {\rm id}|\ \lambda\in F\}$.

It is known that the complexification of  $sl_n(\mathbb R)$ is equal to $sl_n(\mathbb C)$, that is, $sl_2(\mathbb R)$ is  absolutely-simple. As an example of a not absolutely simple Lie algebra, one can consider the special complex Lie algebra  $\mathfrak g=sl_2(\mathbb C)$ and consider $\mathfrak g$ as an algebra over $\mathbb R$ (there is a non-scalar centrum map in $\mathrm{Cent}(\mathfrak g)$, see Example 2). In this case, the complexification $\overline{\mathfrak g}=\mathfrak g_{\mathbb R}\otimes \mathbb{C}$ is a semisimple Lie algebra (\cite{Goto}).

The following theorem is a corollary of Theorem 3 (wich was originally proved in \cite{AndJan} for $F=\mathbb R$). 

{\bf Theorem 5.} Let $\mathfrak g$ be a simple finite-dimensional Lie algebra over a field $F$ of characteristic 0. Then any Lie bialgebra structure on $\mathfrak g$ is either triangular, factorizable, or almost-factorizable. 

{\bf Proof.} 
Let $\delta:\mathfrak g\rightarrow \mathfrak g\otimes \mathfrak g$ be a Lie bialgebra structure on $\mathfrak g$. By the first Whitehead’s Lemma, any Lie bialgebra structure $\delta$ on $\mathfrak g$ is cobaundary, that is, there exists $r\in \mathfrak g\wedge\mathfrak g$ such that $\delta=\delta_r$.

 Let $R_r$ be the corresponding to $r$ map defined by \eqref{TenzOp} and $\theta$ be the corresponding bilinear map defined by \eqref{p41}. By Corollary 4, there is a map $\nu\in\mathrm{Cent}(\mathfrak g)$ such that
$$
\theta(x,y)=\nu([x,y])
$$
for all $x,y\in \mathfrak g$. Since the kernel of $\nu$ is an ideal of $\mathfrak g$, either $\nu=0$ or $\nu$ is non-degenerate. If $\nu=0$, then by Theorem 3, $\delta_r$ is a triangular Lie bialgebra. 

Let $\nu\neq 0$. If $-\nu=\mu^2$ for some $\mu\in \mathrm{Cent}(\mathfrak g)$, then  $(\mathfrak g, \delta)$ is a factorizable Lie bialgebra by Theorem 1. Otherwise, consider the algebraic closure $\bar{F}$ of the field $F$. Let $\mathfrak g_{\bar{F}}=\mathfrak g\otimes _F \bar{F}$ and $\bar{\delta}_r=(\delta_{r})_{\bar{F}}:\mathfrak g_{\bar{F}}\rightarrow \mathfrak g_{\bar{F}}\otimes \mathfrak g_{\bar{F}}$ be the Lie bialgebra structure on $\mathfrak g_{\bar{F}}$ induced by $\delta_r$. Obviously, $\bar{\delta}_r$ is a coboundary Lie bialgebra. Let $\bar{\theta}$ and $\bar{\nu}$ be the corresponding maps defined by $\eqref{p41}$ and Corollary 4. 

Since $\mathfrak g$ is a simple finite-dimensional Lie algebra,  $\mathfrak g_{\bar{F}}$ is a finite-dimensional semisimple Lie algebra over an algebraically closed field $\overline{F}$. Let $\mathfrak g_{\overline{F}}=\bigoplus_{k=1}^s I_k$, where $I_j$ are simple ideals of $\mathfrak g_{\overline{F}}$. Obviously, for any $j\in \{1,\ldots, s\}$, $\overline{\nu}(I_j)=I_j$ and the restriction $\overline{\nu}_j=\overline{\nu}|_{I_j}$ is a centrum map defined on $I_j$.  By Schur's lemma, $\overline{\nu}_j=\lambda_j {\rm id}$ for some $\lambda_j\in \overline{F}$, $\lambda_j\neq 0$ . Since $\overline{F}$ is algebraically closed, there are $\alpha_i\in \overline{F}$ such that $\alpha_i^2=-\lambda_i$. Then the map $\mu$ acting as
$$
\mu(x_1+x_2+\ldots+x_s)=\alpha_1 x_1+\alpha_2 x_2+\ldots+\alpha_s x_s,
$$
$x_j\in I_j$ is a centrum map  satisfying $\mu^2=-\overline{\nu}$. Thus, $(\mathfrak g_{\overline F},\overline{\delta}_r)$ is a factorizable Lie bialgebra by Theorem 3, and $(\mathfrak g, \delta)$ is an almost-factorizable Lie bialgebra.

{\bf Proposition 6 (\cite{Chloup}).} Let $\mathfrak g$ be a finite-dimensional simple, not absolutely simple, real Lie algebra. Then any Lie bialgebra on $\mathfrak g$ is either triangular or factorizable.

{\bf Proof.} Consider $\mathrm{Cent}(\mathfrak g)$. Since the kernel of any centrum map is an ideal in $\mathfrak g$, any nonzero element from $\mathrm{Cent}(\mathfrak g)$ is invertible. Therefore, $\mathrm{Cent}(\mathfrak g)$ is a field containing $\mathbb R$ as a proper subfield. Since $\mathrm{Cent}(\mathfrak g)\subset \mathrm{End}(\mathfrak g)$ and $\mathfrak g$ is a finite-dimensional algebra, $\mathrm{Cent}(\mathfrak g)$ is a finite extension of $\mathbb R$. Thus, $\mathrm{Cent}(\mathfrak g)$ is isomorphic to $\mathbb C$, and for any $\nu\in \mathrm{Cent}(\mathfrak g)$, there exists $\mu\in \mathrm{Cent}(\mathfrak g)$ such that $\mu^2=\nu$. By Theorem 3, any Lie bialgebra on $\mathfrak g$ is either triangular, or factorizable.  The proposition is proved.

{\bf Example 6.} Let $\mathfrak g=sl_2(\mathrm R)$ with the standard basis $x,h,y$ (see Example 2). Consider $r=x\otimes y-y\otimes x$. From \cite{Gomez}, we know that up to the action of $\mathrm{Aut}(\mathfrak g)$, there are three non-isomorphic Lie bialgebra structures given by $r_1=h\wedge x$ (triangular), $r_2=x\wedge y$ (factorizable), and $r_3=h\wedge (x+y)$ (almost-factorizable). 

Let $\omega$ be the bilinear form on $\mathfrak g$ defined by $\omega(x,y)=1$, $\omega(h,h)=2$. For the illustration of the presented theory, consider corresponding to $r_i$ linear maps $R_i$ \eqref{TenzOp} and bilinear maps $\theta_i$ \eqref{p41}.

We have that $R_1(x)=0$, $R_1(y)=-h$ and $R_1(h)=2x$. Then 
$$\theta_1(h,y)=[R_1(h),R_1(y)]-R_1([R_1(h),y]+[h,R_1(y)])=4x-R(2h)=0.
$$

Similarly, $\theta_1(x,h)=\theta_1(x,y)=0$. Thus, $(\mathfrak g,\delta_{r_1})$ is a triangular Lie bialgebra by Theorem 3.

Consider $r_2$. We have that $R_2(x)=-x$, $R_2(y)=y$, $R_2(h)=0$. Then
$$
\theta_2(x,y)=[R_2(x),R_2(y)]-R_2([R_2(x),y]+[x,R_2(y)])=-h-R_2(h-h)=-h=-[x,y],
$$
$$
\theta_2(x,h)=[R_2(x),R_2(h)]-R_2([R_2(x),h]+[x,R_2(h)])=-R_2(2x)=2x=-[x,h].
$$
Similarly, $\theta_2(y,h)=-[y,h]$. Thus, $\theta_2(a,b)=-[a,b]$ for all $a,b\in \mathfrak g$ and $(\mathfrak g, \delta_2)$ is a factorizable Lie bialgebra by Theorem 3.

Consider $r_3$. In this case, $R_3(x)=R_3(y)=-h$, $R_3(h)=2(x+y)$. We have
\begin{gather*}
\theta_3(x,y)=[R_3(x),R_3(y)]-R_3([R_3(x),y]+[x,R_3(y)])=-R_3(2y+2x)=4h=4[x,y],\\
\theta_3(x,h)=-4x+4y-R_3(2h)=-8x=4[x,h],\\
\theta_3(y,h)=-4x+4y-R_3(-2h)=8y=4[y,h].
\end{gather*}
Thus, $\theta_3(a,b)=4[a,b]$ for all $a,b\in \mathfrak g$. Thus, $(\mathfrak g,\delta_{r_3})$ is a Lie bialgebra by Theorem 2. Since $\mathfrak g$ is absolutely simple, $\mathrm{Cent}(\mathfrak g)$ is isomorphic to $\mathbb R$. By Theorem 3, $(\mathfrak g,\delta_{r_3})$ can't be a triangular or a factorizable Lie bialgebra. By Theorem 5, $(\mathfrak g,\delta_{r_3})$ is an almost-factorizable Lie bialgebra.   

{\bf Example 7} Let $\mathfrak g=sl_2(\mathbb R)\oplus sl_2(\mathbb R)$ be the simple Lie algebra from Example 2. Then $\mathfrak g$ is a 6-dimensional simple, not absolutely simple, real Lie algebra with the basis $\{x,h,y,ix,ih,iy\}$, where $\{x,h,y\}$ is the standard basis of $sl_2(\mathbb R)$ (\cite{Goto}). Let $\omega$ be a bilinear map defined as
$$
1= \omega(x,y)=\frac{1}{2}\omega(h,h)=-\omega(ix,iy)=-\frac{1}{2}\omega(ih,ih).
$$

Simple computations show that $(\mathfrak g,\omega)$ is a quadratic Lie algebra. Consider an element $r=h\wedge (x+y)-(i\otimes i)(h\wedge (x+y))$, the corresponding map $R$ defined as
$$
R(x)=R(y)=-h,\ R(h)=2(x+y),\ R(ix)=R(iy)=-ih,\ R(ih)=2i(x+y).\ 
$$
Note that  we have the following equation: $R(ia)=iR(a)$ for all $a\in \mathfrak g$. From the previous example, we obtain that the corresponding bilinear map $\theta$ (defined by \eqref{p41}) satisfies
$$
\theta(a,b)=4[a,b]=-(2\varphi_i)^2([a,b])
$$
for all $a,b\in \mathfrak g$, where $\varphi_i$ is the centrum map from Example 2 defined as $\varphi_i(u+iv)=iu-v$ for all $u,v\in sl_2(\mathbb R)$. Thus, $(\mathfrak g, \delta_r)$ is a factorizable Lie bialgebra. In order to find the corresponding solution of CYBE, consider the map
$B=\frac{1}{2}(R-2\varphi_i)$. We have that
\begin{gather*}
B(x)=-\frac{1}{2}h-ix=-iB(ix),\ B(y)=-\frac{1}{2}h-iy=-iB(iy),\\
B(h)=x+y-ih=-iB(ih).
\end{gather*}
By Statement 2, $B$ is a Rota-Baxter operator of weight $2\varphi_i$. Since $R_r$ is skew-symmetric, $B$ satisfies $B+B^*+2\varphi_i=0$. By Theorem 1, the corresponding to $B$ tensor 
$$
r_B=-(\mathrm{id}\otimes \varphi_i+\varphi_i\otimes \mathrm{id})(\frac{1}{2}h\otimes h+x\otimes y+y\otimes x)
+\frac{1}{2}(\mathrm{id}\otimes\mathrm{id}-\varphi_i\otimes \varphi_i)(h\wedge (x+y)).
$$
is a solution of CYBE on $\mathfrak g$ with $\mathfrak g$-invariant symmetric part.

The following proposition was originally proved in \cite{Chloup} for finite-dimensional simple real Lie algebras.

{\bf Proposition 7.} Let $\mathfrak g$ be a simple finite-dimensional Lie algebra over a field of characteristic 0, $r\in \mathfrak g\wedge\mathfrak g$, and $R_r\in \mathrm{End}(\mathfrak g)$ be the corresponding to $r$ map defined by \eqref{TenzOp}. Then \\
(1) $(\mathfrak g,\delta_r)$ is a triangular Lie bialgebra if and only if the classical double $D(\mathfrak g)$ is not a semisimple Lie algebra. In this case, $D(\mathfrak g)=A\oplus D$ (a direct sum of vector spaces), where $D$ is the radical of $D(\mathfrak g)$ such that $D^2=0$.\\
(2) $(\mathfrak g,\delta_r)$ is a factorizable Lie bialgebra if and only if $D(A)$ is a semisimple not simple Lie bialgebra isomorphic to $\mathfrak g\oplus \mathfrak g$.\\
(3) $(\mathfrak g,\delta_r)$ is an almost-factorizable Lie bialgebra if and only if $D(\mathfrak g)$ is a simple, not absolutely simple, Lie algebra.

{\bf Proof.} Suppose that $(\mathfrak g,\delta_r)$ is a Lie bialgebra. Let $D(\mathfrak g)$ be the classical double of $(\mathfrak g,\delta_r)$. By Proposition 3, $D(\mathfrak g)=D_{R_r,0}(\mathfrak g)$. Let $I$ be a proper ideal in $D(\mathfrak g)$. By Proposition 1, $\mathrm{dim}(I)=\mathrm{dim}(\mathfrak g)$. Recall that $D(\mathfrak g)$ is a quadratic Lie algebra with the nondegenerate bilinear form $Q$ defined as $Q(x+f,y+g)=f(y)+g(x)$ for all $x,y\in \mathfrak g$, $f,g\in \mathfrak g^*$. Let $I^{\perp}$ be the orthogonal complement of $I$ in $D(\mathfrak g)$ with respect to the form $Q$. Since $(D(\mathfrak g),Q)$ is a quadratic Lie algebra and $I$ is a proper ideal, $I^{\perp}$ is a proper ideal in $D(\mathfrak g)$. Therefore, $\mathrm{dim}(I^{\perp})=\mathrm{dim}(\mathfrak g)$. The intersection $I_1=I\cap I^{\perp}$ is also a proper ideal in $D(\mathfrak g)$. Therefore, we have two options: $I=I^{\perp}$ or $I\cap I^{\perp}=0$. 

If $I=I^{\perp}$, then $Q(I^2,D(\mathfrak g))\subset Q(I,[I,D(\mathfrak g)])=0$. Thus, $I^2=0$ and $(\mathfrak g, \delta_r)$ is a triangular Lie bialgebra (see, for example, \cite{GME2}). 

If $I\cap I^{\perp}=0$, then $I^2\neq 0$, $(I^{\perp})^2\neq 0$, and $D(\mathfrak g)=I\oplus I^{\perp}$. Since any proper ideal of $I$ or $I^{\perp}$ is an ideal in $D(\mathfrak g)$, $I$ and $I^{\perp}$ are simple Lie algebras. Therefore, $D(\mathfrak g)$ is a sum of two simple ideals, and $(\mathfrak g,\delta_r)$ is a factorizable Lie bialgebra (see, for example, \cite{GME1}). Now we can use Corollary 5 to complete the proof.

We want to finish with the following result 

{\bf Theorem 6.} Let $\mathfrak g$ be a finite-dimensional semisimple Lie algebra over an algebraically closed field of characteristic 0, $\omega:\mathfrak g\times\mathfrak g\rightarrow F$ is a nondegenerate invariant biliner form on $\mathfrak g$, and $\delta:\mathfrak g\rightarrow \mathfrak g\wedge \mathfrak g$ is a linear map. Then $(\mathfrak g,\delta)$ is a Lie bialgebra if and only if $\delta=\delta_r$ for some $r\in \mathfrak g\otimes \mathfrak g$ such that the corresponding linear map $R_r$ defined by \eqref{TenzOp} is a Rota-Baxter operator of a weight $\mu\in \mathrm{Cent}(\mathfrak g)$ satisfying
$$
R_r+R_r^*+\mu=0.
$$

{\bf Proof.} Suppose that $(\mathfrak g, \delta)$ is a Lie bialgebra. By the first Whitehead’s
Lemma, any Lie bialgebra structure $\delta$ on $\mathfrak g$ is cobaundary, that is, there exists $q\in \mathfrak g\wedge \mathfrak g$ such that $\delta=\delta_q$. Let $R_q$ and $\theta:\mathfrak g\times \mathfrak g\rightarrow \mathfrak g$ be the maps defined by \eqref{TenzOp} and \eqref{p41}. By Theorem 2, $\theta(x,y)=\mu([x,y])$ for all $x,y\in \mathfrak g$, where $\mu\in\mathrm{Cent}(\mathfrak g)$. Using a similar technique as in Theorem 5, we can find a centrum map $\nu\in\mathrm{Cent}(\mathfrak g)$ such that $\nu^2=-\mu$. By Theorem 3, $(\mathfrak g,\delta_r)$ is either triangular or quasitriangular Lie bialgebra. By Theorem 1, $\delta_q=\delta_r$, where $r\in \mathfrak g\otimes\mathfrak g$ and the corresponding map $R_r$ is a Rota-Baxter operator of a weight $\mu\in \mathrm{Cent}(\mathfrak g)$ satisfying \eqref{l11}.

The inverse statement follows from Theorem 1. The theorem is proved.

\section*{Acknowledgments}
The research is supported by Russian Science Foundation (project 23-71-10005, https://rscf.ru/project/23-71-10005/).

\end{document}